\documentstyle[epsf]{article}
\newcommand{\be}{\begin{equation}}
\newcommand{\ee}{\end{equation}}
\newcommand{\bea}{\begin{eqnarray}}
\newcommand{\eea}{\end{eqnarray}}
\newcommand{\beas}{\begin{eqnarray*}}
\newcommand{\eeas}{\end{eqnarray*}}
\newcommand{\bdm}{\begin{displaymath}}
\newcommand{\edm}{\end{displaymath}}
\newcommand{\ba}{\begin{array}}
\newcommand{\ea}{\end{array}}
\newcommand{\bi}{\begin{itemize}}
\newcommand{\ei}{\end{itemize}}
\newcommand{\ben}{\begin{enumerate}}
\newcommand{\een}{\end{enumerate}}
\newcommand{\bc}{\begin{center}}
\newcommand{\ec}{\end{center}}
\newcommand{\bfl}{\begin{flushleft}}
\newcommand{\efl}{\end{flushleft}}
\newcommand{\bfr}{\begin{flushright}}
\newcommand{\efr}{\end{flushright}}
\newcommand{\bd}{\begin{description}}
\newcommand{\ed}{\end{description}}
\newcommand{\bq}{\begin{quote}}
\newcommand{\eq}{\end{quote}}
\newcommand{\bfg}{\begin{figure}}
\newcommand{\efg}{\end{figure}}
\newcommand{\bt}{\begin{table}}
\newcommand{\et}{\end{table}}
\newcommand{\btb}{\begin{tabular}}
\newcommand{\etb}{\end{tabular}}
\newcommand{\btg}{\begin{tabbing}}
\newcommand{\etg}{\end{tabbing}}

\newcommand{\qslash}
           {\mbox{$ q \hspace{-1.1ex} \mbox{/} \hspace{-0.05ex} $}}

\begin{document}
{\small\bf \noindent MZ-TH//97-25\\
\noindent  q-alg/9707029\\[2cm]}
\begin{center}
{\Large {\bf On the Hopf algebra structure of perturbative 
quantum field theories\\[1cm]}}
{\large Dirk Kreimer\footnote{Heisenberg Fellow\\ 
dirk.kreimer@uni-mainz.de\\ 
http://dipmza.physik.uni-mainz.de/$\sim$kreimer/homepage.html}
}\\[7mm]
{\em  Dept.~of Physics}\\
{\em  Univ.~of Mainz}\\
{\em  Staudingerweg 7}\\
{\em  55099 Mainz}\\
{\em  Germany\\[2cm]}
\end{center}
\begin{abstract}
We show that the process of renormalization encapsules
a Hopf algebra structure in a natural manner. This sheds light on
the recently proposed connection between knots and renormalization theory.
\end{abstract}


\section{Introduction}
In this paper, we want to address algebraic properties of renormalization.
Let us first motivate such an attempt. There is a natural grading on Feynman diagrams
by their loop number. But as far as renormalization is concerned there is a 
further grading which is natural to consider. This is
the grading by the number of subdivergences in the graph.
Let us agree to say that an overall divergent graph without 
any subdivergences is of grade one. Then, an overall divergent graph 
with $n-1$ subdivergences is of grade $n$. It is an analytic
expression in which we can localize $n-1$ subgraphs which are
themselves overall divergent. We will see that we can build up
such a graph from $n$ graphs, each of them having no subdivergences
itself.



Renormalization is an operation on graphs which compensates all
(sub-)divergences by help of the forest formula. 
For a given graph $\Gamma$ this formula constructs 
an expression $\bar{\Gamma}$, a quantity which has all its subdivergences
renormalized. From this quantity one obtains a counterterm expression
$Z[\bar{\Gamma}]$ such that $\bar{\Gamma}+Z[\bar{\Gamma}]$ is a finite
quantity. The counterterm is at most a polynomial expression
in external
parameters like masses and momenta and thus is a local
operator.
Its specific form depends on the
chosen renormalization scheme, and on the regularization as well.
All this has become standard textbook knowledge. 


Nevertheless,
renormalization is widely regarded as an {\em ad hoc} procedure.
It is indisputible though that it encapsules some of the most typical
properties of QFT.  In spite of this it is often
considered as a feature which exhibits our present lack of understanding
of QFT proper, and merely regarded as a conceptual asset of our theory.


In this paper we want to demonstrate that renormalization
is not such an {\em ad hoc} procedure, but on the contrary
arises in a very natural
manner from the properties of (quasi-) Hopf algebras.
Our point of departure will be a closer look 
at the forest formula.
This formula is an operation  which is concerned with subgraphs
$\gamma_i\in \Gamma$,
and with replacing the graph $\Gamma$ by a sum
of expressions of the form
$R[Z[\gamma_i]\Gamma/\gamma_i]$, where $\gamma_i$ is a proper subgraph
and $\Gamma/\gamma_i$ is an expression where in the
big graph $\Gamma$ we shrink $\gamma_i$ to a point.
So, renormalization is concerned with replacing a graph $\Gamma$ by
products of graphs, and summation over all possible subgraphs.
Both, $\bar{\Gamma}$ and the corresponding counterterm
$Z[\bar{\Gamma}]$ have this form.
Their difference is finite. 


This finiteness means that we can establish an equivalence
relation
\bea
A\sim B \Leftrightarrow \lim_{\hbar \to 0}(\lim_{\epsilon \to 0} [A-B])=0.
\eea
Thus, we consider expressions which have the same pole part in the regularization
parameter $\epsilon$, but different finite terms, as equivalent.
We are interested in the renormalization procedure
{\em per se}, and not in the finite differences
between various renormalization schemes.



Now consider an equation taken from the theory of Hopf algebras.
\bea
m[(S\otimes id)\Delta[X]]\sim E\circ\bar{e}[X]=0.
\eea
It is our aim to  bring the whole theory of renormalization
in a form such that 
\begin{itemize}
\item
there exists a counit $\bar{e}$ which annihilates
Feynman graphs,
\item the coproduct $\Delta$ generates all the terms necessary to compensate
the subdivergences,
\item in the set ${\cal A}$ furnished by all Feynman diagrams we identify
graphs without subdivergences as the primitive elements in the Hopf algebra,
\item the antipode $S[R[X]]$ coincides with the counterterm $Z_{[X]}$
for a given graph $X$,
\item renormalization schemes correspond to
maps $R:{\cal A}\to {\cal A}$ which fulfil $R[X]\sim X$,
\item $m[(S\otimes id)\Delta[X]]$ is the finite term ($\sim 0$) which corresponds to the
renormalized Feynman graph.
\end{itemize}


To give an idea what is to come, let us consider the two-loop graph $\Gamma^{[2]}(q)$. 
It belongs to the class $((x)x)$ in Fig.(3).
It is the first member of this
class, a one-loop vertex subdivergence $\Gamma^{[1]}(q)$
nested in a one-loop vertex correction, in a theory with a log-divergent
vertex, say.
So it is overall divergent and has a single nested
subdivergence.
Its renormalization reads
\bea
Z_{\Gamma^{[1]}} & = & -R[\Gamma^{[1]}(q)],\\
\overline{\Gamma^{[2]}}(q) & = & \Gamma^{[2]}(q)-Z[\Gamma^{[1]}(q)]\Gamma^{[1]}(q),\\
Z_{\Gamma^{[2]}} & = & -R[\overline{\Gamma^{[2]}(q)}]=
-R[\Gamma^{[2]}(q)]+R[Z[\Gamma^{[1]}(q)]\Gamma^{[1]}(q)],\\
\overline{\Gamma^{[2]}}(q)+Z_{\Gamma^{[2]}}
 & = & \Gamma^{[2]}(q)-Z[\Gamma^{[1]}(q)]\Gamma^{[1]}(q)
-R[\Gamma^{[2]}(q)]+R[Z[\Gamma^{[1]}(q)]\Gamma^{[1]}(q)]  \sim  0.
\eea
The last line in this equation we will show to be generated as
\be
m[(S\otimes id)\Delta[\Gamma^{[2]}(q)]]\sim 0,
\ee
where the coproduct $\Delta$ generates terms $R[\Gamma^{[2]}(q)]\otimes e$, 
$R[\Gamma^{[1]}(q)]\otimes \Gamma^{[1]}(q)$ and
$e \otimes \Gamma^{[2]}(q)$,  the antipode delivers all the correct signs,
and the multiplication map $m$ multiplies everything together to give the renormalized
finite expression.


The antipode of $\Gamma^{[2]}$ then equals $Z_{\Gamma^{[2]}}$ in this example.
This remains true for any given Feynman diagram:
its antipode equals its counterterm.


We will soon see that our 
approach collects Feynman diagrams in classes: all Feynman diagrams
which have a similar forest structure realized by similar subgraphs belong to the same
class. 
Elements in our Hopf algebra correspond to classes of Feynman diagrams,
and any set of Feynman diagrams can be ordered by these classes.



We will see that the primitive elements of the Hopf algebra
are the Feynman diagrams without subdivergences.
This Hopf algebra describes renormalization. It is neither
cocommutative nor coassociative, but most likely these failures are
determined by ${\cal R}$-matrices and associators, and unit-constraints.\footnote{The
reader should not confuse the renormalization map $R$
with a potential ${\cal R}$-matrix for our Hopf algebra. 
The study of the latter is postponed to future work.}
In fact, we will show that there are
renormalization schemes $R$ in which our Hopf algebra is proper.
The coproduct is coassociative and there exists a well-defined counit.
It remains non-cocommutative though.
We refer to \cite{Kassel} for the definition of these notions. 




It is precisely the presence of divergences in the Feynman graphs
which prohibits the classical limit $\hbar\to 0$.
Without the nontrivial renormalization property, the resulting Hopf
algebra were trivial.


The rest of the paper is organized as follows.


We will first study the forest formula and
explain its connection to 
Hopf algebra structures. 
Then we discuss
the Hopf algebra of renormalization
in simple toy models. Finally we treat 
realistic quantum field theories.
In a concluding section, we comment how these findings relate to the recent 
association of knot and number theory to Feynman diagrams.
\section{Feynman diagrams as a realization}
Let us remind ourselves of the forest formula \cite{Collins,Zimm}.
The forest formula organizes the compensation of divergences in Feynman 
graphs.
Such divergences can appear when all the loop-momenta in a graph
tend to infinity jointly. This constitutes the overall divergence
of a graph. But divergences can also appear in sectors where some
loop momenta are kept fixed, while others tend to infinity.
These are proper subdivergences. 
Divergent sectors in a Feynman graph are determined by a powercounting
on all its possible subgraphs \cite{Collins}.


Each divergent subgraph constitutes a forest.
The compensation of divergences
from subgraphs is achieved by constructing a 
counterterm for each divergent subgraph. After we compensated each
subgraph by a counterterm, we are left with the overall divergence
in the graph. Subtracting this final overall divergence determines
the overall counterterm for the Feynman graph.


All this is expressed by the forest formula which we write as
\be
 Z(\Gamma)=-R(\Gamma)-\sum_{\gamma\in \Gamma}R(Z(\gamma)\;\Gamma/\gamma). 
\ee
Here, we use the following notation.
$Z$ determines the overall counterterm of a Feynman graph $\Gamma$.
We sum over all graphs $\gamma$ which are proper subgraphs of $\Gamma$.
If $\gamma$ has no overall divergence, we set $Z(\gamma)=0$.
Otherwise,
$Z(\gamma)$ replaces $\gamma$ in $\Gamma$.
The resulting expression multiplies the graph $\Gamma/\gamma$.
This graph is achieved by reducing $\gamma$ in $\Gamma$ to a point.


$R(\Gamma)$ is a renormalization map.
In the MS scheme, for example, it simply maps to
the proper pole part of $\Gamma$.
For graphs $\gamma_i$ without subdivergences, we have $Z(\gamma_i)=-R(\gamma_i)$.


The above formula applies iteratively until each forest is expressed in terms
of such graphs.
On subgraphs $\gamma$ containing several disconnected components,
e.g.~$\gamma=\gamma_1\cup\gamma_2$, we have $Z(\gamma)=Z(\gamma_1)Z(\gamma_2)$.


The reader can easily check the example in the previous section
against this formula.
We now give a more complicated example.
We consider Fig.(1).
\begin{figure}
\epsfxsize=10cm
\epsfbox{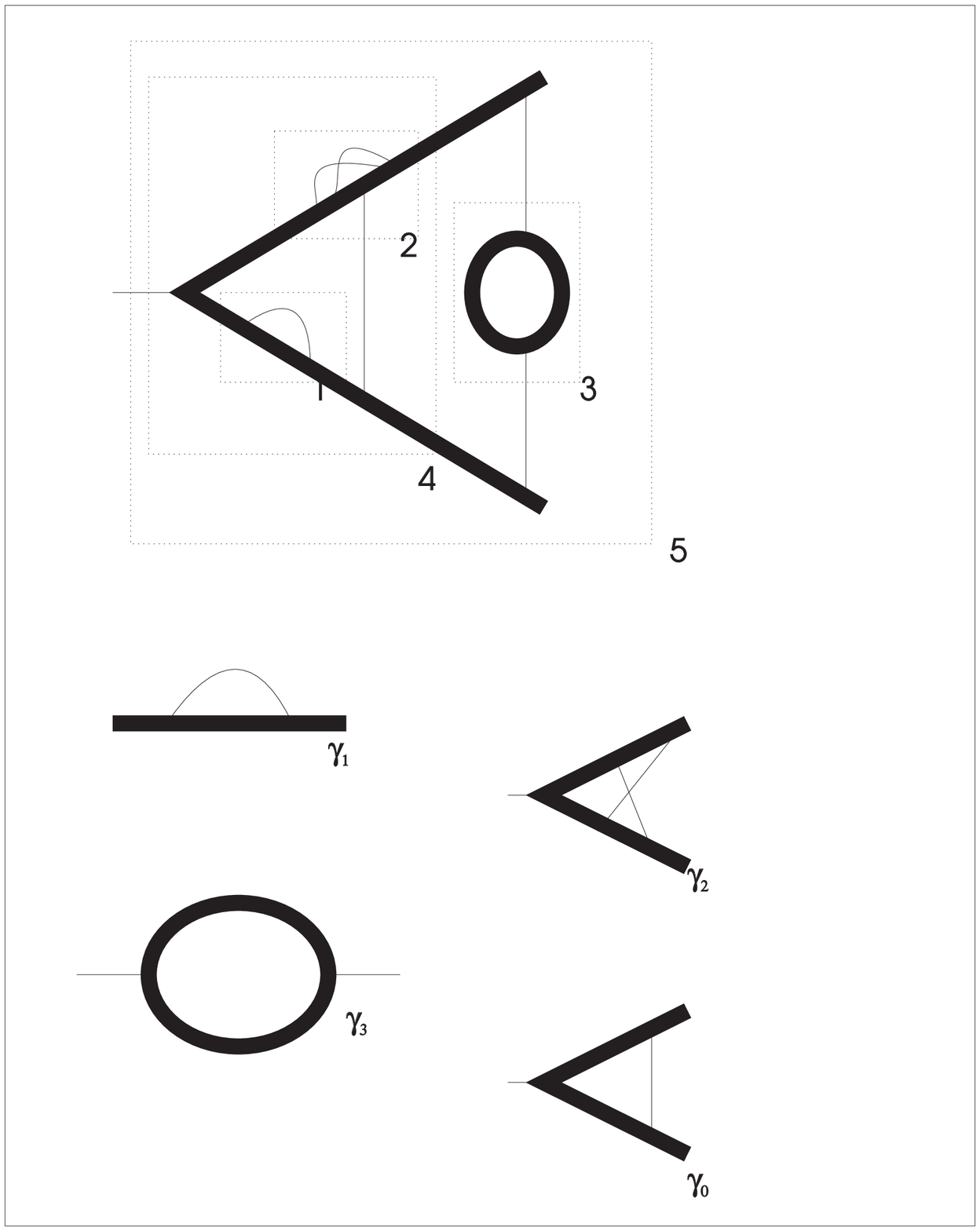}
\caption{A Feynman graph $\Gamma$
with subdivergences. There are various subgraphs
$\gamma_i$. Ultimately, we can consider $\Gamma$ as an expression
in various Feynman graphs without subdivergences,
which are nested into each other and in this way built up the given graph.
The dashed boxes indicate five forests, each containing (sub)-divergences.
The box indicated by the number 5 contains the whole graph $\Gamma$,
while boxes $1,\ldots,4$ contain subdivergences $\gamma_1,\ldots,\gamma_4$.
Note that box number 4 contains the graph $\gamma_4$ which is itself built up
from graphs $\gamma_1,\gamma_2$ and $\gamma_0$.}
\end{figure}
We find
\beas
Z(\Gamma) & = & -R\left(\gamma_5)
-R(Z(\gamma_3)\Gamma/\gamma_3
+Z(\gamma_4)\Gamma/\gamma_4
+Z(\gamma_2)\Gamma/\gamma_2
+Z(\gamma_1)\Gamma/\gamma_1\right.\\
 & & 
+Z(\gamma_3)Z(\gamma_4)\Gamma/(\gamma_3\cup\gamma_4)\\
 & & 
+Z(\gamma_3)Z(\gamma_2)\Gamma/(\gamma_3\cup\gamma_2)\\
 & & 
+Z(\gamma_3)Z(\gamma_1)\Gamma/(\gamma_3\cup\gamma_1)\\
 & & 
+Z(\gamma_3)Z(\gamma_2)Z(\gamma_1)\Gamma/(\gamma_3\cup\gamma_2\cup\gamma_1)\\
 & & 
\left.
+Z(\gamma_2)Z(\gamma_1)\Gamma/(\gamma_2\cup\gamma_1)\right)\\
Z(\gamma_4) & = & -R(\gamma_4)-R(Z(\gamma_1)\gamma_4/\gamma_1+
Z(\gamma_2)\gamma_4/\gamma_2+Z(\gamma_1)Z(\gamma_2)\gamma_4/
(\gamma_1\cup\gamma_2)),
\eeas
with
\be
Z(\gamma_i)=-R(\gamma_i),\;i\in\;\{1,2,3\}.
\ee
We see that we resolve the forests until we express everything in terms of 
Feynman graphs 
without subdivergences. 


It is a standard result that forests never overlap. They appear always
nested inside each other
($\gamma_1$ is inside $\gamma_4$) or disjoint
($\gamma_4$ is disjoint to $\gamma_3$). It is one of the achievements
of renormalization theory that overlapping divergences are resolved in terms of 
nested and disjoint ones. 


For such non-overlapping forests,
we can map our Feynman graph to words with brackets as follows.
Words consists of letters.
Letters are provided by Feynman graphs without
subdivergences. Brackets in such words are derived from
the forest structure and indicate the way subdivergences are entangled
in the Feynman graphs.
To derive the bracket structure from the forest
structure we use the fact that
a Feynman diagram drawn in the plane together with all its forests
is invariant under diffeomorphisms of the plane.
First of all, we draw the diagrams such that
all propagators and vertices which constitute a subdivergence (and only those)
are completely drawn in the interior of the corresponding forest,
for all subdivergences.
Now we use the diffemorphism invariance to deform such a diagram and
its forests in the plane.
This means that we can draw any Feynman diagram with $r$
subdivergences in a manner such that
\begin{itemize}
\item
its forests are rectangular boxes which never overlap,
\item
the sides of these boxes are parallel to the $x$ and $y$ axis,
\item
the outermost box which contains the whole diagram and all
other boxes is the unit square
with corners $(0,0),(0,1),(1,0),(1,1)$,
\item every other box has a height $h$ with $0.9<h<0.95$, say,
and width $1/(r+1)$.
\end{itemize}
Now we delete the horizontal lines of each box apart from a small
reminder (much smaller than the width) close to the corners.
Thus, what remains from a box looks like a pair of brackets.
No such bracket is empty.


We can now map the picture so obtained to words with brackets.
Any pair of bracket might contain other brackets coming from
forests inside. Consider a given pair of brackets.
If we shrink all brackets (forests) inside it to a point,
the remainder is a Feynman graph $\gamma_i$ without subdivergences.
We write the corresponding letter next to the right closing bracket of the
pair of brackets under consideration. The rest of what is contained 
in the pair under consideration is drawn to the left of this letter.


Note that disjoint forests and thus configurations
in disjoint pairs of brackets commute in this construction,
as the notion of being to the left or right is interchangeable
by diffeomorphisms of the plane. 
On the other hand, note that boxes which are inside other boxes
will so remain, as the notion of being inside
is invariant under diffeomorphisms.



Let us have a closer look at our example. What we just
said means that we  can actually
arrange all divergent sectors in a way such that 
forests  are placed almost at the same horizontal level, 
as in Fig.(2).
\begin{figure}
\epsfxsize=10cm
\epsfbox{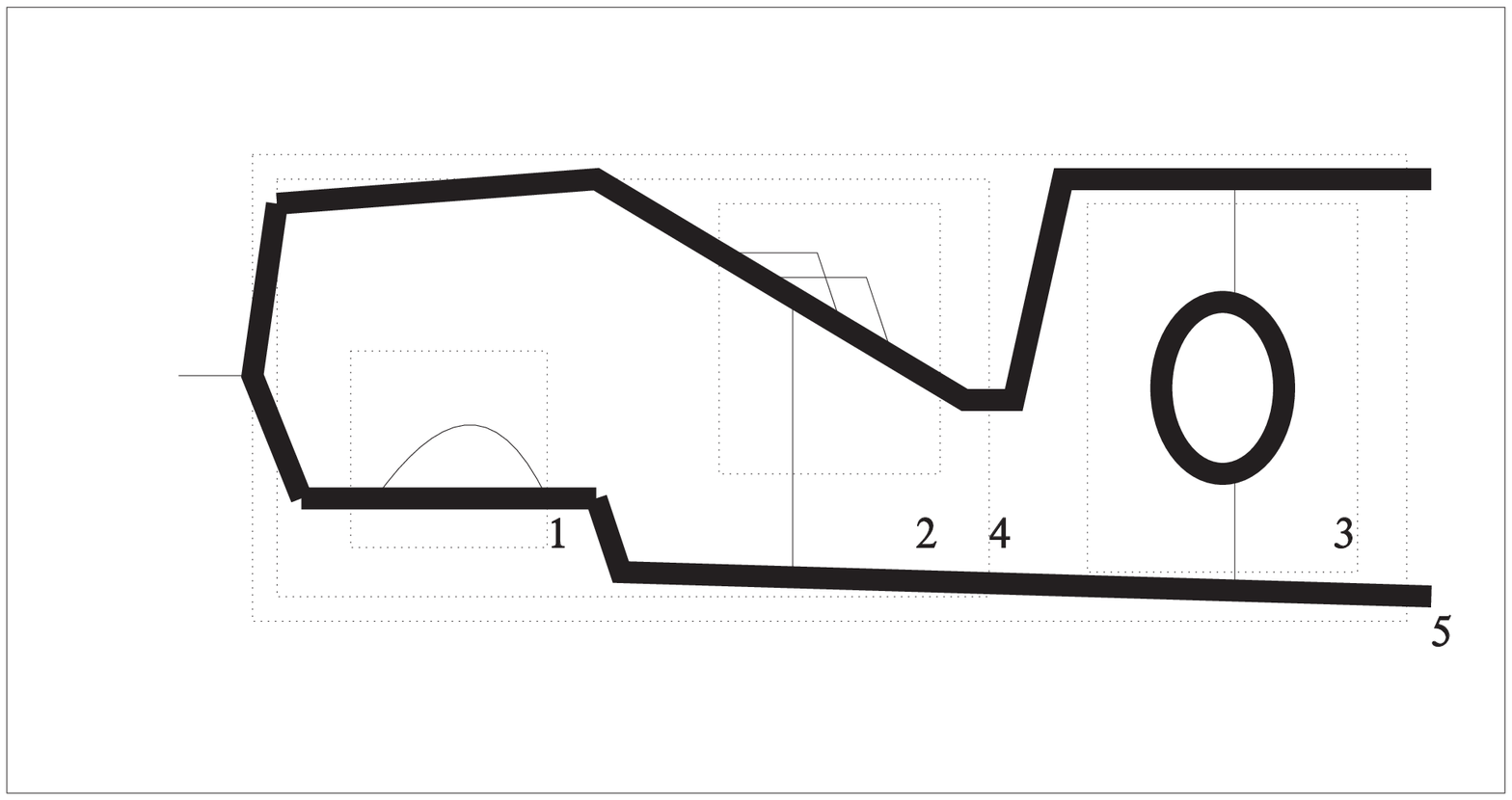}
\caption{We redraw the Feynman graph of the previous figure.
This Feynman graph belongs to the class spanned by the word
$(((\gamma_1)(\gamma_2)\gamma_0)(\gamma_3)\gamma_0)$,
by the procedure explained in the text.
The bracket configuration is obtained from the forest configuration
above by ignoring the horizontal lines in the dashed boxes.
Note that the forests $4,5$ both contain the graph $\gamma_0$ as the
outermost letter on the rhs. Accordingly, if we shrink all subdivergences
in these forests to a point, it is the graph $\gamma_0$ which remains.}
\end{figure}


Now we can safely forget about the horizontal lines 
and interprete the left vertical line of a forest as an
opening bracket, and the right vertical line as a closing one.
Then, we can map the above Feynman diagram to
a word with parenthesis.
It consists of letters which are
Feynman diagrams free of subdivergences:
\be
\Gamma \to (((\gamma_1)(\gamma_2)\gamma_0)(\gamma_3)\gamma_0),
\ee 
while the subgraph $\gamma_4\in ((\gamma_1)(\gamma_2)\gamma_0)$.
Note that 
\be
\gamma_0=\Gamma/(\gamma_3\cup\gamma_4)
\ee
and also
\be
\gamma_0=\gamma_4/(\gamma_1\cup\gamma_2)
\ee 
in our example.




The brackets indicate the forest structure of the diagram.
In this mapping from Feynman diagrams to words with
parenthesis we loose some information. This concerns the location
as to which propagator or vertex of
a graph $\gamma_j$
a graph $\gamma_i$ is attached  as a subdivergence.
There can be several different attachments which generate the same
forest structure.
Naturally, our map to words with parenthesis will
not depend on the choice where to locate the attachment.


But as we have to sum over all possible Feynman diagrams anyhow,
this ambiguity is resolved automatically. Any Feynman diagram lies
in the preimage of some word with parenthesis.


Thus, any Feynman diagram
belongs to a class given by some word with brackets.
A complete class is the set of all Feynman diagrams which
have contributions in a given word with bracket configuration.
Note that due to the very construction of our map disjoint forest configurations
are commutative: $((\gamma_1)(\gamma_2)\gamma_0)=((\gamma_2)(\gamma_1)\gamma_0)$.



In the next section, we will see that there is a natural Hopf algebra
structure on these words, and that the forest formula is naturally related
to the antipode of this Hopf algebra.
\section{The Hopf Algebra}
Consider an alphabet of letters $x_i$. We allow for an
infinite (though countable) set of such letters.
Each Feynman diagram without subdivergences will provide such a letter.
Usually, we have an infinite number of such letters provided by 
a renormalizable theory. A superrenormalizable theory would actually provide
only an alphabet consisting of a finite number of letters.
Our approach is applicable as long as the theory is renormalizable
in the modern sense: a finite number of subtractions is sufficient
at a given loop order, though the number of subtractions (the number of
monomials in the Lagrangian which need to be renormalized)
is allowed to depend on the number of loops considered.


As we have seen, any Feynman diagram belongs to a class of Feynman diagrams.
All elements in this class share that they have the same forest structure.
We describe now such classes by words. Such words consists of letters and round
brackets. The brackets are the remnants of the forest structure.
We agree to write the letter which corresponds to the overall divergence of
a given expression always on the rhs, so that 
the example of the introduction
would correspond to $\Gamma^{[2]}(q)=((x_1)x_1)$ (if $x_1$
is the letter for $\Gamma^{[1]}$) and {\em not} to
$(x_1(x_1))$.
Thus we will have classes of Feynman diagrams
represented by  the notion of parenthesized words (PWs) on letters $x_i$.
Here are some examples
\be
(),(x_i),((x_i)x_j),(x_i)(x_j),((x_i)(x_j)x_k),\ldots
\ee
From what we have said so far it is clear 
that each PW consists of a string of letters and balanced brackets
such 
that each letter in a PW has {\em
one and only one closing bracket on its rhs}, but
can have more than one opening bracket on its lhs.
The empty word $e:=()$ is allowed as well and will act as
the unit in the algebra to be defined in a moment.


A PW whose leftmost opening bracket is matched by its 
rightmost closing bracket is called an irreducible PW (iPW).
Otherwise it is reducible and a product of iPWs.
Here are some irreducible words
\be
(),(x_i),((x_i)x_j),((x_i)(x_j)x_k),(((x_i)x_j)x_k),\ldots
\ee
and here is a reducible word
\be
((x_i)x_j)(x_k).
\ee
For an arbitrary PW we call the number of letters in it
its length $k$.
As each letter corresponds to a Feynman diagram
with some given loop number, and without subdivergences, 
each such letter  is assumed to be of order
\be
x_i\sim \hbar^{N_i}/\epsilon,
\ee
where $N_i$ is its loop number. 
Thus, $N_i$ is a positive integer $(\geq 1)$,
which we assume to be known for each $x_i$,
and $\hbar,\epsilon$
are small real positive parameters. 
This is input taken from perturbative quantum field theory.
Feynman graphs without subdivergences all provide a first
order pole in a regularization parameter $\epsilon$,
and are of order $\hbar^{N_i}$, if they contain $N_i$ loops.


Thus, each PW delivers a Laurent series in $\epsilon$.
Its degree equals the length of the PW.
Further, the order or numbers of loops $N$
of a PW of length $k$ is the sum of the orders in $\hbar$ of all its
letters,
\be
N:=\sum_{j:=1}^k N_{i_j}.
\ee
Let us also introduce the depth of a iPW.
It is the maximum number of nested divergences in a iPW,
counting the final overall divergence as well.
For example, the following words have depths 1,2,2,3 from the left to
the right: $(x),((x)x),((x)(x)x),(((x)x)x)$.
They correspond to a graph $x$ without subdivergence, one with a single
subdivergence, one with two disjoint subdivergences, and one with
a subdivergence which has a nested subdivergence itself.
For the first, second and fourth example the depth and the length actually
agree. We call such iPWs strictly nested.



Let ${\cal A}$ be the set of all PWs.
We regard it as a ${\bf Q}$-vectorspace. 
We furnish this vectorspace with a product $m$ and a unit $e\equiv()$,
$eX=Xe=X$,
by introducing
\bea
m: & & {\cal A}\otimes {\cal A}\to{\cal A}\nonumber\\
m[ X\otimes Y ] & := & X \;Y \equiv Y\;X,\\
E: & & {\bf Q}\to {\cal A}\nonumber\\
E(q) & := & q e, 
\eea
$\forall$ PWs $X$, $Y$ and rational numbers $q$
as linear maps.


A word about notation is appropriate.
To avoid the proliferation of round brackets, we write
capital letters $X,Y,\ldots$ for PWs. One should not forget that
they stand for expressions which have brackets indicating the forest structure.
Small letters like $x,y,\ldots$ indicate single letters from an alphabet
provided by Feynman diagrams without subdivergences.
Their brackets are always given explicitly. 
An example: if $X$ stands for $((x)x)$ then $(Xy)$ is the well-defined iPW $(((x)x)y)$.





Let us introduce an equivalence relation at this point.
As pointed out in the introduction, we are mainly concerned with the
divergent parts of graphs. Any Feynman graph delivers a Laurent series
in a regularization parameter $\epsilon$, and we consider two expressions
as equivalent if these Laurent series have identical pole parts.
Hence we call two expressions $X,Y$ equivalent,
$X\sim Y$, if and only if 
\be
\lim_{\hbar\to 0,\epsilon\to 0} [X-Y]=0.
\ee
We extend this equivalence relation to ${\cal A}^{\otimes n}$ such that,
for example,
\be
[X\otimes Y]\sim [U\otimes V]\Leftrightarrow
\lim_{\hbar\to 0,\epsilon\to 0}[X-U]=0
\;\mbox{and} \;
\lim_{\hbar\to 0,\epsilon\to 0}[Y-V]=0 \;\mbox{and}\;
\lim_{\hbar\to 0,\epsilon\to 0}[XY-UV]=0,
\ee
at the same time.
The limits have to be taken in the order as indicated.


An endomorphism $R$ 
from ${\cal A}$ $\to$ ${\cal A}$
is called a renormalization map if and only if
\be
[R[X]-X]\sim 0.
\ee
Accordingly, $R[X]$ has the same pole terms in $\epsilon$ as $X$.
Various different renormalization schemes correspond to the choice of
different maps $R$.
With this definition we introduce a renormalization scheme as a map
which leaves the divergences of a primitive Feynman graph
(one without subdivergences) unchanged. Thus it generates a 
proper counterterm which can annihilate the divergence of the
primitive graph. General graphs are then to be handled by the
combinatorics of the Hopf algebra.


Note that in general 
\be
[R[X]-X]\otimes Y\not\sim 0,
\ee
as Y provides pole terms which forbid the limit $\epsilon\to 0$.


So far we have furnished the set ${\cal A}$ with an algebra structure.
Next, we furnish ${\cal A}$ with a Hopf algebra structure \cite{Kassel}.
We define a counit by
\bea
\bar{e}: & & {\cal A}\to {\bf Q}\nonumber\\
\bar{e}[e] & := & 1,\\
\bar{e}[X] & := & 0,\;\forall X\not=e\;\in {\cal A}.
\eea
Before we can check that this is a counit we have to define a coproduct
$\Delta$. A few words of caution are appropriate. 
In the literature, one often speaks about
Hopf algebras and quasi-Hopf algebras.
While the former have a proper counit and
a coassociative coproduct, the latter provide
such structures only up to isomorphisms.
The coproduct we will define in a moment will be proven to be
coassociative only for some special cases of $R$.
Only then it is a proper coproduct, and ${\cal A}$ a proper
Hopf algebra. A similar remark applies to the counit.
For general $R$, it fails to be a counit.
There are hints that all these failures are quasi.
This means that we hope to be able to define left and right unit constraints
and an non-trivial associator in the case of general $R$. For the case studied here
we are content to restrict ourselves to the case of a {\em coassociative} $R$.
By this we mean that $R$ fulfils a constraint given below, which makes
${\cal A}$ into a proper Hopf algebra. 
It is very gratifying that there are renormalization schemes which 
fulfil such a constraint. 
For such renormalization schemes, we have a proper Hopf algebra
at hand, as we will show.


Let us continue with the introduction of a coproduct. 
The coproduct is a sum of terms $X_i\otimes Y_i$. Roughly speaking, the
terms on the lhs $X_i$ are considered to be subdivergent graphs. Somehow
we must stop producing counterterms if $X_i$ happens to be the letter e,
that is, when there is no divergent graph as the factor on the lhs.
This is done by a mapping $P_L$ defined as follows.
\bea
P_L: & & {\cal A}\otimes{\cal A}\to{\cal A}\otimes{\cal A}\nonumber\\
P_L & := & (id-E\circ \bar{e})\otimes id
\eea
Note that $P_L[e\otimes X]=0,$ $ \forall$ PWs $X$, and that
$P_L^2=P_L$. 


Also, we want our coproduct to be able to act on the
subdivergences and to plug the result of its action back into the
graph. This is achieved with the help of the following endomorphism
 $B_{(x_i)}$ of ${\cal A}\otimes{\cal A}$
which depends itself on a single letter $x_i$:
\bea
B_{(x_i)}: & & {\cal A}\otimes{\cal A}
\to{\cal A}\otimes{\cal A}\nonumber\\ 
B_{(x_i)}[X\otimes Y] & := & X\otimes(Yx_i).
\eea
An example might be in order:
\be
B_{(x_2)}[R[(x_1)]\otimes (x_3)]=R[(x_1)]\otimes ((x_3)x_2).
\ee
The letter $x_i$ which parametrizes $B$ corresponds to the Feynman graph
which we obtain when we shrink all subdivergences to a point.


With the help of the maps $P_L$ and $B$ we are now prepared to define
the coproduct:
\bea
\Delta: & & {\cal A}\to{\cal A}\otimes {\cal A}\nonumber\\
\Delta[e] & := & e\otimes e\label{cop0}\\
\Delta[(x_i)] & := & R[(x_i)]\otimes e+e\otimes (x_i)\label{cop1}\\
\Delta[X\; Y] & := & \Delta[X]\Delta[Y],\\
\Delta[(Xx_i)] & := & R[(Xx_i)]\otimes e + e \otimes
(Xx_i)+B_{(x_i)}[P_L[\Delta[(X)]]],\label{copm}\\
\Delta(R[X]) & := & \Delta[X].\label{de4}
\eea
Note that the above definition is complete
as any iPW has the form $(Xx_i)$ for some not necessarily irreducible
PW $X$ and for
some letter $x_i$.
The presence of $P_L$ guarantees that we defined the coproduct
on PWs of length $k$ in terms of the coproduct of words of length $<k$.
Actually, the above definition is slightly redundant.
Eq.(\ref{cop1}) follows from Eqs.(\ref{cop0},\ref{copm}).
Also, there are various other ways to write this coproduct.
For example, adopting Sweedler's notation (\cite{Kassel}),
we find
\bea
\Delta[(Xx)] & = & R[(Xx)]\otimes e + \sum_{X}R[X^\prime]\otimes (X^{\prime\prime}x)\\
 & = & R[(Xx)]\otimes e + B_{(x)}\left[\sum_{X}R[X^\prime]\otimes X^{\prime\prime}\right].\label{cop}
\eea


The reader should use these definitions 
to confirm the following examples, which are given
for the case of
a Feynman graph with one nested subdivergence
in Eq.(\ref{ex1}),  for two disjoint subdivergences
nested in a further graph in Eq.(\ref{ex2}), 
and finally for a subdivergence nested in another graph
nested in yet another graph in Eq.(\ref{ex3}).
Fig.(3) contains examples for these cases.
\begin{figure}
\epsfxsize=12cm
\epsfbox{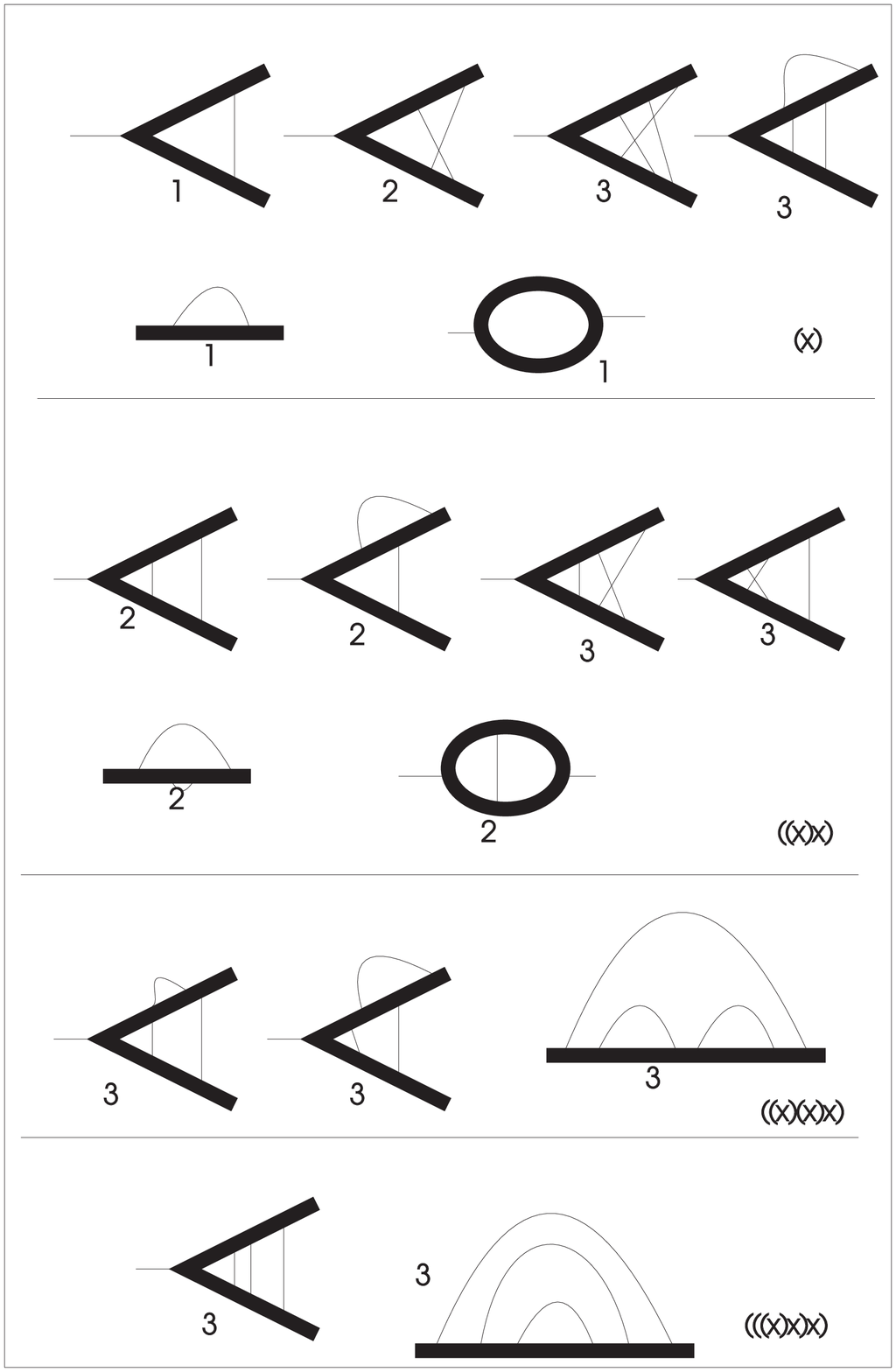}
\caption{Some Feynman graphs belonging to the classes $(x)$, $(x)x)$,
$((x)(x)x)$ and $(((x)x)x)$. We indicate the order in $\hbar$.}
\end{figure}


\bea
\Delta[((x_i)x_j)] & = & R[((x_i)x_j)]\otimes e + e \otimes ((x_i)x_j)
+R[(x_i)]\otimes (x_j),\label{ex1}\\
\Delta[((x_i)(x_j)x_k)] & = & R[((x_i)(x_j)x_k)]\otimes e 
+ e \otimes ((x_i)(x_j)x_k)\nonumber\\
 & &
+R[(x_i)]\otimes ((x_j)x_k)
+R[(x_j)]\otimes ((x_i)x_k)\nonumber\\
 & & 
+R[(x_i)]R[(x_j)]\otimes (x_k),\label{ex2}\\
\Delta[(((x_i)x_j)x_k)] & = & R[(((x_i)x_j)x_k)]\otimes e 
+ e \otimes (((x_i)x_j)x_k)\nonumber\\
 & & 
+R[(x_i)]\otimes ((x_j)x_k)
+R[(x_i)x_j)]\otimes (x_k).\label{ex3}
\eea
Having defined a coproduct, it is natural to ask if it is
cocommutative. We will also have to check
its coassociativity later on.
Obviously, $\Delta$ is not cocommutative.
Flipping the rhs of Eq.(\ref{ex1}), we get
\bea
\Delta[((x_i)x_j)]-\Delta^{op}[((x_i)x_j)] & \sim & 
R[(x_i)]\otimes (x_j)-(x_j)\otimes R[(x_i)]\;\not\sim  0.\label{exop}
\eea


It is now time to check the counit.
The most general element of ${\cal A}$ has the form
$X=\prod_i R[X_i]\prod_j Y_j$ for iPW's $X_i,Y_j$.
We note that the coproduct of $X$
always has the
form $\prod_i \prod_jR[X_i]R[Y_j]\otimes e+e\otimes X+\ldots$, 
where the dots represent terms
which neither have a unit $e$ on the lhs or rhs of the tensor product.
These terms are annihilated by $(id\otimes \bar{e})$ as well as by $(\bar{e}\otimes
id)$.
Thus we obtain
\bea
(id\otimes \bar{e})\Delta[X] & = & \prod_i R[X_i]\prod_jR[Y_j]\not\sim X\\
(\bar{e}\otimes id)\Delta[X] & = & \prod_i X_i\prod_j Y_j\not\sim X.
\eea
This indicates a problem.
It shows that $\bar{e}$ is not a proper counit, for
general $R$. There are two possible solutions.
One either abandons the notion of a Hopf algebra, and tries to
establish the apparatus of a quasi Hopf algebra for general 
maps $R$. 
This will be done in future work \cite{km}.


For our present purposes we are more modest. 
We will formulate a further condition
on $R$ such that we can remain in the setup of a Hopf algebra.
We will see that the condition which makes $\bar{e}$ a counit
also guarantees coassociativity.
Clearly, the condition is
\be
R[\prod_i R[X_i]\prod_j Y_j]= R[\prod_i X_i \prod_j Y_j].\label{cond}
\ee
This repairs the counit, as
\be
\prod_i R[X_i]\prod_j Y_j\sim
R[\prod_i R[X_i]\prod_j Y_j]= R[\prod_i X_i \prod_j Y_j]
\sim \prod_i X_i\prod_j Y_j.
\ee
We will see later that nontrivial maps $R$ ($R=id$ would be a trivial choice)
which are in accordance with
Eq.(\ref{cond}) exist.


Next we study coassociativity. We want to show that for renormalization
maps $R$ which fulfil the above condition, $\Delta$ is coassociative.
For general renormalization maps $R$, $\Delta$ is not coassociative,
which one can easily show by considering $(\Delta\otimes id)\Delta
-(id\otimes\Delta)\Delta$ acting on words of lengths $>1$.
But the coassociativity is restored if
$R$ fulfils Eq.(\ref{cond}). 
Let us prove this assertion.


$\Delta$ is obviously coassociative on words of length one for all $R$.
For higher words, we procced by induction.
We use the presentation of the coproduct as given in Eq.(\ref{cop}).
For the induction, we assume that $\Delta$ is coassociative on words of length $n$.
We first prove that it is coassociative on irreducible words
of length $n+1$ and then prove the assertion for general PW's.
A $\sum_X$ in accordance with
 Sweedler's notation is understood in what follows. 
\beas
(id\otimes \Delta)\Delta[(Xx_j)] & = & 
(id\otimes \Delta)B_{(x_j)}[ R[X^\prime]\otimes X^{\prime\prime}]
+R[(Xx_j)]\otimes e\otimes e\\
 & = & 
(id\otimes B_{(x_j)})(id\otimes \Delta)\Delta[X]+ R[X^\prime] \otimes R[(X^{\prime\prime}x_j)]
\otimes e\\
 & & 
+R[(Xx_j)]\otimes e\otimes e\\
 & = & 
(id\otimes B_{(x_j)})(\Delta\otimes id)\Delta[X]+R[X^\prime]\otimes
R[(X^{\prime\prime}x_j)]\otimes e\\
 & & 
+R[(Xx_j)]\otimes e\otimes e\\
 & = &
(id\otimes B_{(x_j)})(\Delta\otimes id)\Delta[X]+\left[B_{(x_j)}( R[X^\prime]\otimes
X^{\prime\prime})\right.\\
 & & 
\left. +R[(Xx_j)]\otimes e\right]\otimes e\\
 & = & 
(\Delta\otimes id)B_{(x_j)}[\Delta[X]]+\Delta[(Xx_j)]\otimes e\\
 & = &
(\Delta\otimes id)(\Delta[(Xx_j)]).
\eeas
We used Eq.(\ref{cop}) several times and the assumption of coassociativity
at depth $n$ in the third line.
In the fourth line we used 
\be
R[X^\prime]R[(X^{\prime\prime}x_j)]=
R[R[X^\prime](X^{\prime\prime}x_j)],
\ee
by  Eq.(\ref{cond}) for $R$
to conclude
\be
R[X^\prime]\otimes R[(X^{\prime\prime}x_j)]\otimes e \sim 
R[X^\prime]\otimes (X^{\prime\prime}x_j)\otimes e.
\ee


It remains to prove that $\Delta$ is coassociative
on arbitrary words of length $n+1$.
\beas
(id\otimes \Delta)\Delta[XY] & = &
(id\otimes \Delta)\Delta[X]\Delta[Y]\\
 & = & 
(id\otimes \Delta)\Delta[X]
(id\otimes \Delta)\Delta[Y]\\
 & = & 
(\Delta\otimes id)\Delta[X]
(\Delta\otimes id)\Delta[Y]\\
 & = &
(\Delta\otimes id)\Delta[XY]. 
\eeas


Thus, the failure of coassociativity can be attributed
to the fact that in general $R$ is not in accordance with Eq.(\ref{cond}). 


Renormalization schemes
which do not fulfil the condition of coassociativity are the MS scheme
or the on-shell scheme. On the other hand, the momentum scheme is a renormalization
scheme which is coassociative, as we will see below.


Note further that if we were to set
\be
\Delta[R[X]]=(id\otimes R)\Delta[X]
\ee
instead of Eq.(\ref{de4}) we would have a coassociative coproduct for arbitrary $R$,
by inspection of the proof above.
Nevertheless there remains the problem with the counit, which still forces us to either maintain
a condition on $R$ or to build up a quasi Hopf algebra.


We note in passing that $\Delta$ remains non-cocommutative even for the trivial case
$R=id$, and even if we symmetrize in letters $x_i,x_j$,
as Eq.(\ref{ex2}) exhibits.


Having defined a proper
counit and a coassociative
coproduct, we have established a bialgebra
structure. What remains is to find an antipode. It turns out that the antipode
is the object which actually achieves the renormalization: it picks
up the terms generated  by the coproduct and combines them in a way which 
coincides with the action of the forest formula.



Explicitly, we find the antipode $S$ as
\bea
S: & & {\cal A}\to {\cal A}\nonumber\\
S[e] & = & e,\\
S[(x_i)] & = & -(x_i),\\
S[R[(x_i)]] & = & -R[(x_i)],\\
S[XY] & = & S[Y]S[X],\\
S[(Xx_i)] & = & -(Xx_i)-m[(id\otimes S)P_2(\Delta[(Xx_i)])],\label{s1}\\
S[R[(Xx_i)]] & = & -R\left[(Xx_i)+m[(S\otimes id)P_2(\Delta[(Xx_i)])]\right],\label{s2}\\
P_2 & := & (id-E\circ\bar{e})\otimes (id - E\circ\bar{e}).\nonumber
\eea
For consistency reasons, for coassociative $R$ fulfilling Eq.(\ref{cond}),
the two equations Eqs.(\ref{s1},\ref{s2}) must be equivalent,
$S[(Xx_i)]\sim S[R[(Xx_i)]]$ or $R[S[(Xx_i)]]= S[R[(Xx_i)]]$.
To show this we use that $\Delta$ is coassociative under these circumstances.
It is then easy to verify that Eq.(\ref{s1}) is indeed equivalent
to Eq.(\ref{s2}). What has to be shown is that
\beas
 & & m[(id\otimes m)(id\otimes id\otimes S)
(P_3)(id\otimes \Delta)(\Delta)]\\
 & & =
m[(m\otimes id)(S\otimes id\otimes id)
(P_3)(\Delta\otimes id)(\Delta)],
\eeas
where
\be
P_3:=(id-E\circ\bar{e})\otimes(id-E\circ\bar{e})\otimes(id-E\circ\bar{e}).
\ee
This now follows by induction.
Using the associativity of $m$,
the coassociativity of $\Delta$  and
for the induction the assumption that the assertion holds for PWs of length $n$,
we can conclude the assertion for PWs of length $n+1$. We can  shift the $m,S$ and 
$\Delta$
operations as needed, as the presence of the $P_3$ operator guarantees that
$S$ acts only on words of length $\leq n$, where we can shift it by assumption.


To start the induction one
easily verifies that the assertion holds for words of
length three by explicit calculation. There is nothing to prove
at length two or length one.


The following examples are instructive and easy to work out 
\bea
S[R[((x_i)x_j)] & = & -R[((x_i)x_j)]+R[R[(x_i)](x_j)]\\
S[R[((x_i)(x_j)x_k)]] & = & -R[((x_i)(x_j)x_k)]+R[R[(x_i)]((x_j)x_k)]\nonumber\\
 &  &  +R[R[(x_j)]((x_i)x_k)]-R[R[(x_i)]R[(x_j)](x_k)].
\eea


Though for non-coassociative renormalization schemes $R$
\be
S[X]-S[R[X]]\not\sim 0,\;X\in {\cal A},
\ee
one can still use this difference to investigate the action of an associator.
We will have much more to say about this in future work \cite{km}.


Returning to coassociative $R$, we  
show as an example that $R[S[((x_i)(x_j)x_k)]] =S[R[((x_i)(x_j)x_k)]]$:
\beas
S[ ((x_i)(x_j)x_k) ] & = & 
-((x_i)(x_j)x_k)- R[(x_i)]S[((x_j)x_k)]-
R[(x_j)]S[((x_i)x_k)]\\
 & & 
-R[(x_i)]R[(x_j)]S[(x_k)]\\
 & = &
-((x_i)(x_j)x_k) - R[(x_i)](-((x_j)x_k)+R[(x_j)](x_k))\\
 & & 
-
R[(x_j)](-((x_j)x_k)+R[(x_i)](x_k))+R[(x_i)]R[(x_j)](x_k)\\
 & = &
-((x_i)(x_j)x_k) +R[(x_i)]((x_j)x_k)\\
 & & 
+R[(x_j)]((x_j)x_k) -R[(x_i)]R[(x_j)](x_k)\\
 & \Rightarrow & R[ S[ ((x_i)(x_j)x_k) ] ]=S[R[((x_i)(x_j)x_k)]].
\eeas


One verifies from its definition that $S$ fullfils
\be
m \circ (S\otimes id) \circ \Delta\sim
m \circ (id \otimes S) \circ \Delta\sim E \circ \bar{e}
\ee
for general $R$, and is thus a proper antipode, for example:
\bea
m[(S\otimes id)\Delta[(Xx)]] & = &  
S[R[(Xx)]]+\sum_X S[R[X^\prime]](X^{\prime\prime}x)
\nonumber\\
 & = & -R[(Xx)]-R\left[\sum_{{X \atop
X^\prime\not= e}}S[R[X^\prime]](X^{\prime\prime}x)\right]
 +\sum_X S[R[X^\prime]](X^{\prime\prime}x)\nonumber\\
 & \sim & 0=E\circ\bar{e}[(Xx)],
\eea
and similarly for $(id \otimes S)$.


As the algebra ${\cal A}$ is commutative, the antipode $S$  fulfils
$S^2=id$. This is clear from a standard argument \cite{Kassel},
but can also be checked directly using Eqs.(\ref{s1},\ref{s2}).


It is an instructive exercise to work out a few more examples and check
the assertions above.



Our coproduct was constructed so that it can also be written as
\be
\Delta[(Xx)]=\sum_{\mbox{\footnotesize all subwords $U$ of $(Xx)$}}R[U]\otimes (Xx)/U,
\ee
where $(Xx)$ and $e$ are considered as subwords as well,
and the antipode is thus given as
\be
S[R[(Xx)]]=-R[(Xx)]-
\sum_{\mbox{\footnotesize all proper subwords $U$ of $(Xx)$}}R[S[R[U]]\; 
(Xx)/U].
\ee


So, by comparison with the forest formula, 
the formulas for the antipode conspire
to give the $Z_{X}$-factor of a Feynman graph $X$
as $Z_{X}=S[R[X]]$ ($=R[S[X]]$ for coassociative $R$)
and the Feynman graph
with renormalized subdivergences $\bar{X}$
has the form that a renormalized Feynman graph -a finite quantity $\bar{X}+Z_{X}$-, 
is obtained as 
\be
m[(S\otimes id)\Delta[X]].
\ee
For $X=((x_i)x_j)$ this expression reads
\be
((x_i)x_j)-R[(x_i)](x_j)-R[((x_i)x_j)]+R[R[(x_i)](x_j)],
\ee
which is evidently finite and corresponds to a Feynman graph
$(x_i)$ nested in another graph $(x_j)$.
We encourage the reader to work out a few more examples
(like the one of Fig.(1))
and prove the equivalence of the antipode with the forest formula of the
previous section in simple examples.

Note further that one can use the split $id\otimes id=P_R+(id\otimes E\circ\bar{e})$
to write
\be
m[(S\otimes id)\Delta[X]]=
m\left[(S\otimes id)P_R[\Delta[X]]\right]+S[R[X]].
\ee
We conclude that $m[(S\otimes id)P_R[\Delta[X]]]$ coincides with the
famous ${\bar R}$ operation in standard renormalization
theory. This will be nicely exemplified in an example below.




\section{A Toy Model}
\subsection{Nested and disjoint divergences}
In this section we want to discuss some toy models. 
The study of these toy models was motivated by
\cite{mills}.
The purpose is to become acquainted
with the Hopf algebra approach. Both, realistic QFT's to be discussed 
in the next section and the models discussed here realize the Hopf algebra
described in the previous section.
We will also discuss some factorization properties which are useful
from a technical viewpoint, and which can be established for proper
QFT's as well.


Consider a set of functions
\be
I_j^{[1]}(c):=\int_0^\infty \frac{x^{-j\epsilon}dx}{x+c}=
B(j\epsilon,1-j\epsilon)c^{-j\epsilon}.
\ee
We assume $c>0$, $1>\!>\epsilon>0$.


These functions will play the role of the letters $x_j$ in the previous
section. Let us call the integrand of 
$I_j^{[1]}(c)$ a '$j$-loop Feynman diagram
without subdivergences, depending on the external parameter $c$'.
Feynman diagrams have to be integrated, which then
corresponds to an 'integrated Feynman diagram' $I_j^{[1]}(c)$.
The integration of Feynman diagrams is usually ill-defined,
which is reflected by the fact that as a Laurent-series
in $\epsilon$, $I_j^{[1]}(c)$ has a pole term $\sim \Gamma(j\epsilon)$.



We are interested in the limit $\epsilon\to 0$.
So we will replace the Feynman-diagram by the diagram plus
a counterterm, whose sole purpose it is to absorb this pole term.
Actually, field theory demands that such counterterms are 
independent of external parameters. This is the mechanism which brings renormalization
into the game.



It is easy to see that 
\be
I_j^{[1]}(c)-I_j^{[1]}(1)=B(j\epsilon,1-j\epsilon)[c^{-j\epsilon}-1],\label{eq1l}
\ee
has no pole term any longer (in $\epsilon$). 
Thus, the second term on the lhs above
is a good $c$-independent  counterterm.
To get non-trivial examples, we define 'Green functions with subdivergences',
and thus non-trivial PW's.


Let us consider the case of a nested subdivergence:
\be
I_{j_1,j_2}^{[2]}(c)
:=
\int_0^\infty dx \frac{x^{-j_2\epsilon}}{x+c}I_{j_1}^{[1]}(x+c).
\ee
The pole terms of this function do not depend
on the presence of the parameter $c$ in its subdivergence. We thus set
$c=0$ in the subdivergence and one indeed verifies that
\be
I_{j_1,j_2}^{[2]}(c)-I_{j_1+j_2}^{[1]}(c)I_{j_1}^{[1]}(1),
\ee
allows for the limit $\epsilon \to 0$. We note that as far as the pole terms
are concerned the Green function  $I_{j_1,j_2}^{[2]}(c)$ factorizes into a
product of $I_{j}^{[1]}$ functions.
Note further that
\be
I_{j_1+j_2}^{[1]}(c)I_{j_1}^{[1]}(1)=
I_{j_1+j_2}^{[1]}(1)I_{j_1}^{[1]}(1)c^{-(j_1+j_2)\epsilon}.
\ee.


But, if we try to compensate the pole terms in this
expression by subtraction of $I_{j_1+j_2}^{[1]}(1)I_{j_1}^{[1]}(1)$,
we find that neither the pole terms cancel, nor do we find independence
of $c$. The limit $\epsilon\to 0$ of
\be
I_{j_1+j_2}^{[1]}(1)I_{j_1}^{[1]}(1)[c^{-(j_1+j_2)\epsilon}-1]
\ee
does not exist.


The cure is readily at hand. We have to compensate the
subdivergence  first.
The counterterm for it is the scale-independent
$I^{[1]}_{j_1}(1)$, by Eq.(\ref{eq1l}). 
Plugging this in amounts to the replacement
\bea
I_{j_1,j_2}^{[2]}(c) & \rightarrow &
I_{j_1,j_2}^{[2]}(c)-I_{j_2}^{[1]}(c)I_{j_1}^{[1]}(1)\nonumber\\
 & \sim & 
I_{j_1+j_2}^{[1]}(c)I_{j_1}^{[1]}(1) -I_{j_2}^{[1]}(c)I_{j_1}^{[1]}(1).\nonumber
\eea
The equation is correct modulo finite terms.
The final expression corresponds to the 'Green-function plus its counterterm
graph'.


It is now easy to see that this combination has pole terms
which are independent of $c$. Thus, evaluating this combination
at $c=1$ gives us another combination with the same pole terms.
Finally, our counterterm reads 
\be
Z_{I^{[2]}_{j_1,j_2}}=
I_{j_1+j_2}^{[1]}(1)I_{j_1}^{[1]}(1) -I_{j_2}^{[1]}(1)I_{j_1}^{[1]}(1).\nonumber
\ee


Let us pause at this stage and make contact to the previous section.
We define $e:=1$.
We identify primitive elements $(x_i)$ with $I_{i}^{[1]}(c)$.
All PW's $X$ will be a function of a scale $c$, $X\equiv X[c]$.
General PW's are obtained by defining 
\bea
(Xx_i)[c]=\int_0^\infty dy \frac{y^{-i\epsilon}}{y+c}X[y+c]
\eea
and the multiplication of PW's $XY$ is simply the
multiplication of functions $X[c]Y[c]$.
Thus, for example, $((x_i)(x_j)x_k)[c]$ is given by
\be
\int_0^\infty dx_3 dx_2 dx_1 \frac{x_3^{-k\varepsilon}}{[x_3+c]}
\frac{x_2^{-j\varepsilon}}{[x_2+x_3+c]}
\frac{x_1^{-i\varepsilon}}{[x_1+x_3+c]}
\ee
and $(((x_i)x_j)x_k)[c]$ by
\be
\int_0^\infty dx_3 dx_2 dx_1 \frac{x_3^{-k\varepsilon}}{[x_3+c]}
\frac{x_2^{-j\varepsilon}}{[x_2+x_3+c]}
\frac{x_1^{-i\varepsilon}}{[x_1+x_2+x_3+c]}.
\ee
Identifying $R$  with the evaluation at
$c=1$,
\be
R[X[c]]=X[1],
\ee 
we can now determine the antipodes of arbitrary 'Green functions' $(Xx)[c]$.
The reader will find it a simple
but instructive exercise to verify that $m[(S\otimes id)\Delta[(Xx)]]$
is indeed finite. This amounts to a proof of renormalizabilty
on the set of PW's defined as above.
We also note that the above choice of $R$ indeed delivers a proper Hopf
algebra, as it fulfils the condition Eq.(\ref{cond}).


An arbitrary Green function is a combination of disjoint
or nested divergences indicated by a $()$-type bracket structure
introduced in the previous section.
The combination $m[(S\otimes id)\Delta[(Xx)]]$ corresponds to the
renormalized Green function associated to the
unrenormalized Green function $(Xx)$. 


Another very instructive toy model is obtained if we choose
\bea
(x_j)[c] & := & I_j^{[1]}(c)  =  \int_c^\infty dy y^{-1-j\epsilon}\\
(Xx_j)[c] & := & \int_c^\infty dy y^{-1-j\epsilon}X[y]\\
(x_j)(x_k)[c] & := & I_j^{[1]}(c)I_k^{[1]}(c)\\
R[(Xx_j)] & := & (Xx_j)[1].
\eea
We leave it as a useful exercise to the reader 
to verify that all the terms contributing to the renormalized
function conspire to a single iterated integral:
\be
m\left[(S\otimes id)\Delta[(Xx)[c]]\right]=\overline{(Xx)}[c]
:=\int_1^c dy y^{-1-j\epsilon} \overline{X}[y],
\ee
which is evidently finite.
This is in a sense the simplest realization of our Hopf algebra.
The fact that it delivers an iterated integral representation will
be commented elsewhere \cite{book}.
Let us close this section by calculating $((x_1)(x_2)x_1)[c]$
as a very instructive and beautiful example.
It corresponds to two disjoint subdivergences, both of
logarithmic nature. Note that they have different powers in $\epsilon$ (different 'loop
order').
\bea
((x_1)(x_2)x_1)[c] & := & 
\int_c^\infty dx_3 x_3^{-1-\epsilon}
\int_{x_3}^\infty dx_2 x_2^{-1-2\epsilon}
\int_{x_3}^\infty dx_1 x_1^{-1-\epsilon},\\
m[(S\otimes id)\Delta[((x_1)(x_2)x_1)[c]\;]] 
 & = & 
((x_1)(x_2)x_1)[c]-R[(x_1)[c]]((x_2)x_1)[c]-R[(x_2)[c]\;]
(x_1)x_1)[c]\nonumber\\
 & & 
+R[(x_1)[c]\;]R[(x_2)[c]\;]
(x_1)[c]
-R[((x_1)(x_2)x_1)[c]\;]\nonumber\\
 & &
+R[R[(x_1)[c]]((x_2)x_1)[c]]
+R[R[(x_2)[c]\;]((x_1)x_1)[c]\;]\nonumber\\
 & &
-R[R[(x_1)[c]\;]R[(x_2)[c]\;]
(x_1)[c]\;]\\
 & = & 
((x_1)(x_2)x_1)[c]-(x_1)[1]((x_2)x_1)[c]-(x_2)[1]
(x_1)x_1)[c]\nonumber\\
 & &
+(x_1)[1](x_2)[1]
(x_1)[c]
-((x_1)(x_2)x_1)[1]+(x_1)[1]((x_2)x_1)[1]\nonumber\\
 & &
+(x_2)[1]((x_1)x_1)[1]
-(x_1)[1](x_2)[1]
(x_1)[1]\\
 &  = &
 \int_c^\infty dx_3 x_3^{-1-\epsilon}
\int_{x_3}^\infty dx_2 x_2^{-1-2\epsilon}
\int_{x_3}^\infty dx_1 x_1^{-1-\epsilon}\label{term1}\\
 & &
-\int_c^\infty dx_3 x_3^{-1-\epsilon}
\int_{x_3}^\infty dx_2 x_2^{-1-2\epsilon}
\int_{1}^\infty dx_1 x_1^{-1-\epsilon}\label{term2}\\
 & &
-\int_c^\infty dx_3 x_3^{-1-\epsilon}
\int_{1}^\infty dx_2 x_2^{-1-2\epsilon}
\int_{x_3}^\infty dx_1 x_1^{-1-\epsilon}\label{term3}\\
 & &
+\int_c^\infty dx_3 x_3^{-1-\epsilon}
\int_{1}^\infty dx_2 x_2^{-1-2\epsilon}
\int_{1}^\infty dx_1 x_1^{-1-\epsilon}\label{term4}\\
 & & -\;\mbox{the four terms above at $c=1$.}\label{termr}
\eea
Note that the four terms above are amongst them free of subdivergences.
They are generated by $m[(S\otimes id)P_R\Delta[((x_1)(x_2)x_1)[c]]]$,
and coincide with the result the $\bar{R}$ operation would have
obtained at this point.

The other four terms will cancel the remaining overall divergences.
All eight terms together are finite. To verify this we add  them to one
iterated integral.
To this end we combine (\ref{term1}) and term (\ref{term2}) to
\be
-\int_c^\infty dx_3 x_3^{-1-\epsilon}
\int_{x_3}^\infty dx_2 x_2^{-1-2\epsilon}
\int_{1}^{x_3} dx_1 x_1^{-1-\epsilon}.\label{term12}\\
\ee
We further decompose (\ref{term3}) as
\bea
-\int_c^\infty dx_3 x_3^{-1-\epsilon}
\int_{1}^\infty dx_2 x_2^{-1-2\epsilon}
\int_{x_3}^\infty dx_1 x_1^{-1-\epsilon} & = &
-\int_c^\infty dx_3 x_3^{-1-\epsilon}
\int_{1}^\infty dx_2 x_2^{-1-2\epsilon}\nonumber\\
 & & \times \left[
\int_{1}^\infty dx_1 x_1^{-1-\epsilon}
-\int_{1}^{x_3} dx_1 x_1^{-1-\epsilon}\right].\label{term3c}
\eea
The first term in (\ref{term3c}) cancels (\ref{term4}) while
the second term combines with (\ref{term12}) to give
\be
\int_c^\infty dx_3 x_3^{-1-\epsilon}
\int_{1}^{x_3} dx_2 x_2^{-1-2\epsilon}
\int_{1}^{x_3} dx_1 x_1^{-1-\epsilon}.\label{term5}\\
\ee
The other four terms (\ref{termr}) give the same as (\ref{term5}), this
time evaluated at $c=1$,
so that we finally have for the sum of the eight terms
\be
-\int_1^c dx_3 x_3^{-1-\epsilon}
\int_{1}^{x_3} dx_2 x_2^{-1-2\epsilon}
\int_{1}^{x_3} dx_1 x_1^{-1-\epsilon},\\
\ee
which is evidently finite if we send
$\epsilon\to 0$. If we were to consider the corresponding situation
for proper Feynman graphs as in Fig.(3) our Hopf algebra would deliver the
renormalization with the same ease by applying $m[(S\otimes id)\Delta]$
to the corresponding word.






\subsection{Overlapping Divergences}
We mentioned in the first section that it is one of the achievements
of renormalization theory that it disentangles overlapping divergences 
in terms of nested and disjoint ones.
In this subsection, we again take recourse to a simplified model
to study the crucial properties, before we turn to quantum field theory
in the next section.


To have a model for overlapping divergences, we have to consider
expressions which are at least linearly divergent.
In analogy to the previous cases
we consider
\bea
J^{[1]}_{j_1}(c):=\int_0^\infty dx \frac{x^{1-j_1\epsilon}}{x+c}
\eea
as well as
\bea
J^{[2]}_{j_1,j_2}(c):=\int_0^\infty dx_1 
\frac{x_1^{1-j_1\epsilon}}{[x_1+c]}
\frac{1}{[x_1+x_2]}
\frac{x_2^{1-j_1\epsilon}}{[x_2+c]}.
\eea
One verifies that the following expression is finite in the limit
$\epsilon\to 0$:
\be
J^{[2]}_{j_1,j_2}(c)
-\int_0^\infty dy 
\frac{y^{1-j_1\epsilon}}{y+c}I^{[1]}_{j_2}(y)
-\int_0^\infty dy 
\frac{y^{1-j_2\epsilon}}{y+c}I^{[1]}_{j_1}(y)
+
J^{[2]}_{j_1,j_2}(0).
\ee
Thus, the divergences of $J^{[2]}_{j_1,j_2}(c)$
are contained in the other three terms above.
The scale independent term $J^{[2]}_{j_1,j_2}(0)$
is irrelevant (and vanishes in dimensional regularization)
while the other two terms can be written as
\be
J^{[1]}_{j_1+j_2}(1)\left[I^{[1]}_{j_1}(1)+I^{[1]}_{j_2}(1)\right]
c^{-(j_1+j_2)\epsilon}.
\ee


Adding the two counterterm graphs $I_{j_1}^{[1]}(1)J_{j_2}^{[1]}(1)
c^{1-j_2\epsilon}+\{j_1\leftrightarrow j_2\}$ for the overlapping subdivergences
in 
$J^{[2]}_{j_1,j_2}(c)$ we find four terms which can be written as
\be
m[(S\otimes id)\Delta[X_o]],
\ee
where 
\be
X_o=((I_{j_1}^{[1]})J_{j_2}^{[1]})[c]+\{j_1\leftrightarrow j_2\}.
\ee
We see that we resolved the overlapping divergence
into two nested terms of the form $((I)J)$.
The results in \cite{habil} demonstrate that this is always
possible.
We will see in the next section how to derive this result
in general.
\section{Quantum Field Theories}
The construction of our Hopf algebra is made in a way as to deliver
the renormalization of a QFT whenever its renormalization
can be achieved by local counterterms who are determined
by powercounting and thus iterate to a forest formula under
renormalization. In this section we will discuss a few technical issues.
They relate to the possibility to decompose Feynman graphs in terms
of primitive elements. 
Such a decomposition we expect to  prove very convenient
in future work.
We employ the fact that eventually
divergences are independent of external parameters, and thus
all counterterms depend on a single scale. 
\subsection{Internal corrections}
Let us first clarify our aim. A Feynman graph $\gamma_i$
which is only a function
of one external momentum $q$ has the form
\be
\sum_i o_i F_i(\epsilon)[q^2]^{-N_i\epsilon}.
\ee
We sum over the various possible spin structures $o_i$.
They can have a polynomial dependence on masses and the external
momentum $q$. The formfactors $F_i$ will be logarithmicly 
overall divergent scalar functions. They can thus be regarded (modulo
UV-finite terms)
as functions of the
single scale $q^2$. An example in place is the fermion propagator in QED say,
which allows for a decomposition
\be
\qslash F_1(\epsilon)[q^2]+m{\bf 1} F_2(\epsilon)[q^2]
\ee
with two spin structures $\qslash,m{\bf 1}$.


Assume there is another graph $\Gamma$ which contains
$\gamma_i$ as a proper subgraph at some propagator.
Let this propagator carry internal momentum $l$ say.
This propagator will furnish a denominator $1/l^2$.
Then, $\Gamma$ can be calculated as 
\be
\sum_i F_i(\epsilon)\times o_i\star\Gamma/\gamma_i \left\{ \frac{1}{l^2}\to
 \frac{1}{[l^2]^{1+N\epsilon}} \right\}.\label{factor}
\ee
Here, the $\star$-product indicates an appropriate insertion of the spin-structures
$o_i$. The propagator obtains a non-integer powers $1+N\epsilon$,
if $\gamma_i$ contained $N$ loops. Here we assumed that we calculated in
dimensional regularization. The reader will find it easy to modify
the argument for other regularization schemes.
 
 
Eq.(\ref{factor}) is the factorization we are after.
We reduced the calculation into a product form, where we allow for non-integer
powers in the denominators of propagators.
We establish a renormalization scheme
in accordance with Eq.(\ref{cond}) if we set 
$R[\gamma_i[q^2]]=\gamma_i[1]$.\footnote{Or, taking into account the presence of
the scale $\mu$ in DR, we set $q^2=\mu^2$.}


But in general, a graph $\gamma_i$ will depend on more than one
external momentum. If it appears as a proper subgraph
in some larger graph $\Gamma$ so that the external momenta of $\gamma_i$
are internal momenta of $\Gamma$ it cannot be nullified.
Thus, the subgraph depends on various
scales and is in general not factorizable.


In such a situation, we proceed as follows.
Assume $\gamma_i$ is a function of $r>1$ external momenta
$l_{i_1},\ldots,l_{i_r}$, which are internal in $\Gamma$:
\be
\gamma_i=\gamma_i(l_{i_1},\ldots,l_{i_r}).
\ee
We make the replacement 
\be
\gamma_i(l_{i_1},\ldots,l_{i_r})\to
[\gamma_i(l_{i_1},\ldots,l_{i_r})
-\gamma_i(0,\ldots,l_{i_j},\ldots,0)]
+\gamma_i(0,\ldots,l_{i_j},\ldots,0).
\ee
The term in brackets on the rhs does not provide a subdivergence any longer,
as the overall divergence of $\gamma_i$ is independent of the
chosen values of the momenta $l_i$ (with the exception that we are not
allowed to nullify all of them), and thus cancels out in the difference.
We thus do not need a factorization property for it, as it does not
generate a subdivergence and thus no forest.


The other term, which still generates a subdivergence,
scales as $(l_{i_j}^2)^{N_i\epsilon}$ and thus
factorizes as above when plugged into $\Gamma$.


The following example demonstrates the idea.
Let $\gamma_1$ be the one-loop vertex correction
in massless Yukawa theory, say, evaluated at zero momentum
transfer. Let $\gamma_{1,q}$ be the same graph with nonvanishing
momentum transfer.
In Fig.(4) we see two graphs $\Gamma_1,\Gamma_2$. 
$\Gamma_1$ belongs
to the class $(((\gamma_1)\gamma_1)\gamma_1)$.
\begin{figure}
\epsfxsize=10cm
\epsfbox{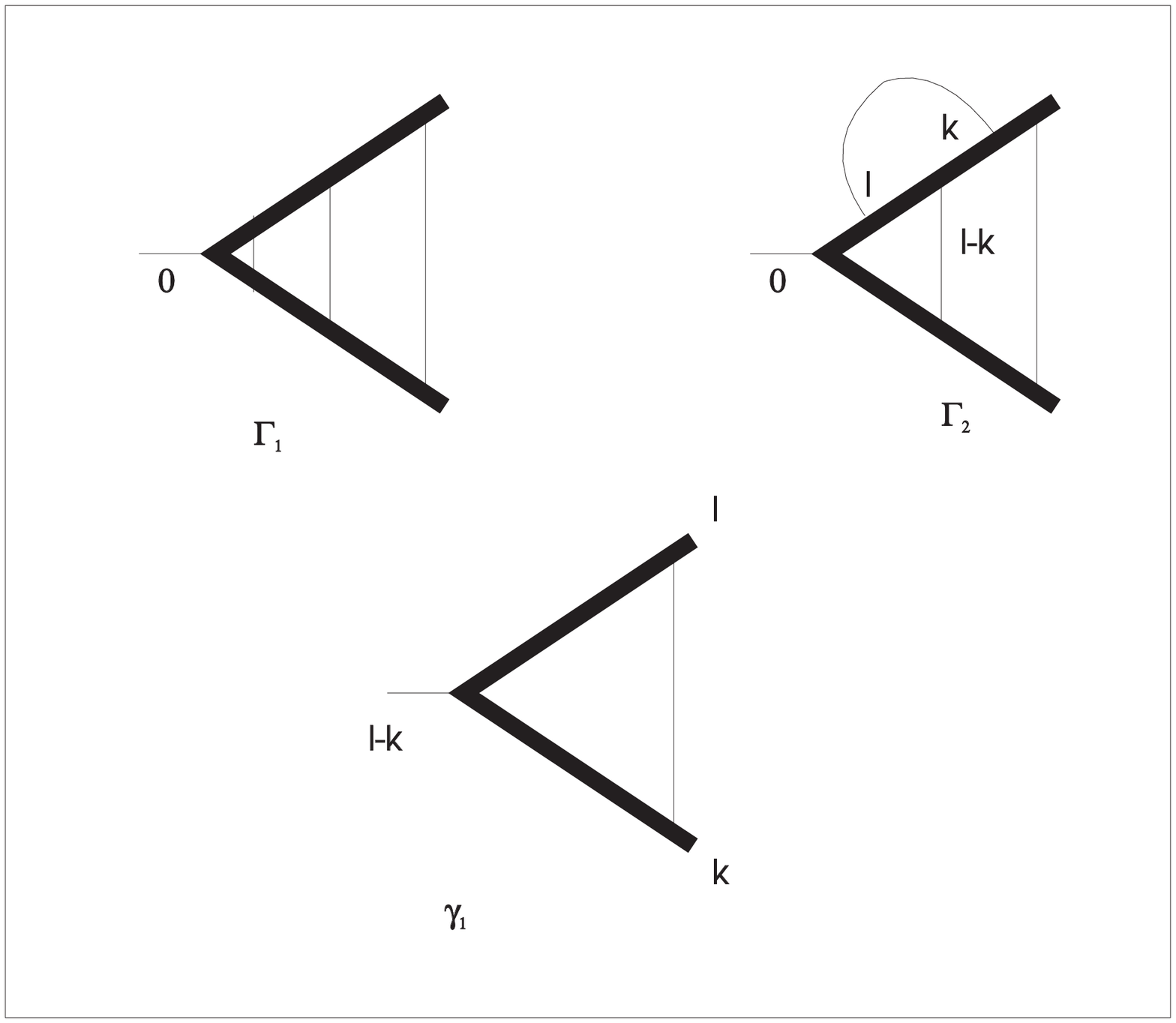}
\caption{Two Feynman graphs. $\Gamma_1$ belongs to the
class $(((\gamma_1)\gamma_1)\gamma_1)$, it is directly factorizable
as we can nullify the external momentum $q$. $\Gamma_2$ splits
into a contribution in this class, plus a contribution
$(X\gamma_1)$, where $X=([\gamma_{1,k-l}-\gamma_1]\gamma_1)$.}
\end{figure}
$\Gamma_2$ splits by the above procedure into two terms:
\be
(([\gamma_{1,q}-\gamma_1]\gamma_1)\gamma_1)+
(((\gamma_1)\gamma_1)\gamma_1).
\ee
We use that $R[\gamma_1]\sim R[\gamma_{1,q}]$
and conclude that 
\be
\Delta[((\gamma_{1,q}-\gamma_1)\gamma_1)]=
R[((\gamma_{1,q}-\gamma_1)\gamma_1)]\otimes e+ e\otimes
((\gamma_{1,q}-\gamma_1)\gamma_1).
\ee
Thus, we can regard $((\gamma_{1,q}-\gamma_1)\gamma_1)$ as a new
primitive element (it has no subdivergences) in its own right. 


It is a general
phenomenon
that internal subdivergences split into factorizable contributions
when expanded in terms of primitive elements.
With this decomposition, 
we have completely reduced the calculation of $Z$-factors (and quantities
derived from them, like anomalous dimensions and $\beta$-functions)
to the calculation of primitive elements. What remains is to classify the
latter, which is a non-trivial task in view of the relations amongs them
\cite{4TR}.


\subsection{Once more: Overlapping Divergences}
In this section, we once more consider overlapping divergences.
We utilize the Schwinger Dyson equation to show that we can resolve
these divergences into nested and disjoint ones in a factorizable manner
as before.
Similar techniques were used in \cite{habil,bdk}.


We start with the Schwinger Dyson equation for a
two-point function.
\bea
P(q) & = & Z\int d^Dk \Gamma(k,q+k)P(q)P(q+k)\nonumber\\
 & = &  
\int d^Dk \Gamma(k,q+k)P(q)P(q+k)\Gamma(k+q,k)\nonumber\\
 & & 
-\int d^Dk \Gamma(k,q+k)P(q)P(q+k)K(k,k+q)P(k)P(k+q)\Gamma(k+q,k).\label{sd}
\eea
In the last line in this equation, we see the presence
of the kernel $K$. It belongs to both subdivergences
and is the source of the overlapping divergences.


Assume now that $P$ is linearly overall divergent.
This allows to improve overall powercounting
by appropriate substitutions like
\be
\Gamma(k,q+k)\to [\Gamma(k,q+k)-\Gamma(k,k)]+\Gamma(k,k)
\ee
and
\be
\Gamma(k,q+k)P(k+q)K(k,k+q)\to [\Gamma(k,q+k)P(k+q)K(k,k+q)-
\Gamma(k,k)P(k)K(k,k)]+\Gamma(k,k)P(k)K(k,k)
\ee
on the lhs and rhs
in the above equations. 
We keep only terms which contain UV-divergences.
This allows to express the overall divergence of $P$ in terms of simpler
functions as demonstrated in Fig.(5). 
\begin{figure}
\epsfxsize=14cm
\epsfbox{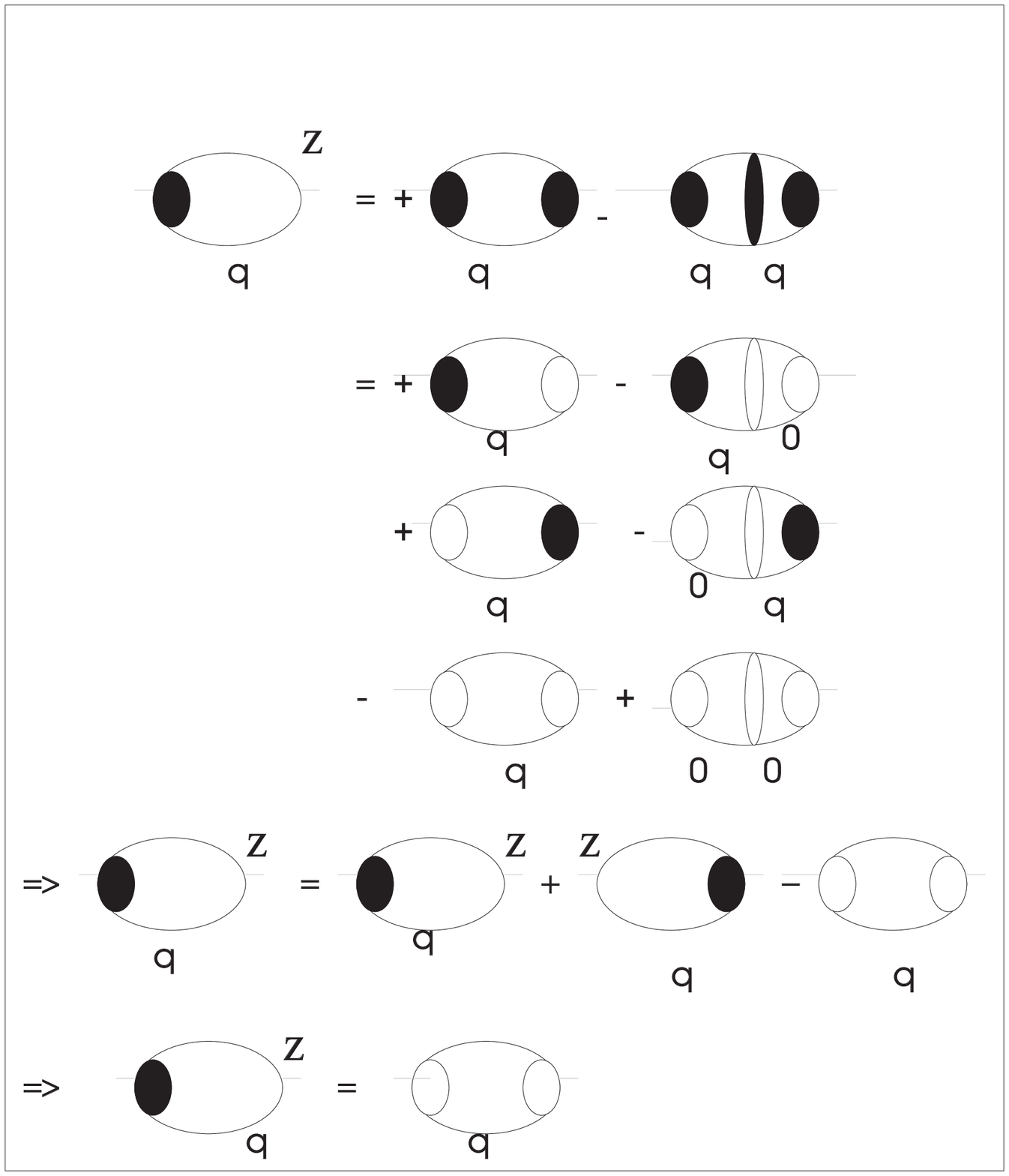}
\caption{Using Eq.(\ref{sd}) back and forwards, we obtain
an easy expression for the overall divergence of an overall
linearly divergent two-point function. This gurantees all necessary
factorization properties and disentangles the overlapping divergences
in accordance with the forest formula.}
\end{figure}


In this figure, we indicate nullified external momenta
by unfilled blobs and a zero instead of $q$ at 
internal propagators.


We see that in the final equation all reference to the kernel $K$ disappeared.
This is the desired result.
This became possible as we could
recombine all terms with the kernel
with other terms to the $Z$-factor of the vertex,
by using the Schwinger Dyson equation for the
vertex. Note that there is one term in Fig.(5) in which all reference
to the external parameter $q$ has disappeared. It corresponds to an expression
without any scale in dimensional regularization and vanishes, as usual.


For quadratic divergent two-point functions, one takes recourse
to the methods of \cite{bdk}, to find similar results.
Using such techniques, one can indeed not only establish that forests never 
overlap (a well-known fact in renormalization theory), but that one can
always resolve subdivergent sectors into factorizable ones
plus contributions in which these sectors become primitive.
\section{Final Remarks}
In this paper we recovered renormalization as a commutative
not cocommutative Hopf algebra.
We establihed that a renormalized finite
quantity is obtained as $m[(S\otimes id)\Delta[X]]$ for a
Feynman diagram $X$ decomposed into a PW.
The overall counterterm $Z_X$ is obtained as the antipode
of $R[X]$. The Hopf algebra is proper for certain renormalization
schemes, distinguished by Eq.(\ref{cond}). 
We expect it to be a braided quasi-Hopf algebra
for general renormalization schemes. Work along these lines
is in progress \cite{km}. 


We expect the algebraic structures described here to be at the heart
of the recently discovered connection between knots, number theory
and Feynman diagrams, especially in the light of the well-known
connection between multiple zeta values and the associator of a
quasi bi-algebra \cite{Kassel}.
Feynman diagrams  evaluate favourably in Euler/Zagier sums
\cite{habil,db,plb,pisa,bdk,bgk,bk15,BBB,db2}. For a definition
of these sums, see for example \cite{db,Zag}. 
In view of the results in \cite{dbn} we are optimistic to relate
this property ultimately to the non-(co-)associativity structures observed in
${\cal A}$.


Let us emphasize a few more interesting aspects.
\begin{itemize}
\item
We organized Feynman diagrams in classes corresponding to their
forest structure. The interesting properties of these classes
come from the fact that nested and disjoint subdivergences
renormalize in a different manner.
The first example appears when one considers elements
\be
(((x)x)x)-((x)(x)x)\not\sim 0.
\ee
The non-triviality of this difference lies at the heart of all
the interesting structures of the Hopf algebra ${\cal A}$. 
Preliminary results indicate that the associator to be established
for general $R$ is determined by differences as the above.
\item
The simplest iPWs are delivered by purely nested configurations
$((\ldots(x_1,)x_1)\ldots )x_1)$
where we only allow for lowest order letters (one-loop graphs) $x_1$.
These words are known to deliver rational antipodes \cite{habil,Bob}.
\item
It is precisely the difference between nested and disjoint subdivergences above
which stops us from defining a well-defined product on graphs.
\item
One can naturally map our PWs to more common
looking bracket configurations, in the spirit of
\bea
(((x)x)x)\to (x\otimes x)\otimes x\\
((x)(x)x)\to x\otimes (x\otimes x).
\eea
This more tensor-categorial way of thinking will be addressed in the
future.
\item
From empirical investigations \cite{db,plb,pisa,bgk,bk15,db2}, we conjecture a deeper role
of positive knots in the understanding of Euler/Zagier sums.
We expect an understanding of ${\cal A}$ to be mandatory in this respect.
\item
It is interesting to note that 
there is a very natural realization of ${\cal A}$
in terms of iterated integrals. We will report consequences elsewhere \cite{book}.
\item
Recently, it was suggested that counterterms in Feynman diagrams establish a weight system.
Subdivergences obscured and restricted the applicability of these results
\cite{4TR,bk4}. The bialgebraic results on
the structure of subdivergences obtained in this paper
will help clarify the situation.
\item
The results reported here fit nicely into the patterns observed in
\cite{sub}. A systematic investigation is ongoing.
\item
For general $R$, $R[S[X]]\not= S[R[X]]$. Nevertheless, one
can show that the two expressions
remain equivalent for general $R$ as long as the depth of $X$ is $\leq$ 2.
\item
David Broadhurst and the author succeeded to turn the Hopf algebra described
here into a computer program which facilitates renormalization.
It will be made available soon \cite{bkhopf}.
\item
The series of different possible forest structures (for irreducible words)
starts as 1,1,2,4,9,20,51,121,
321,826,2186,5789,16114,42449,$\ldots$ \cite{sloane}.
It follows an easy recursion. We will report on its properties
elsewhere \cite{book}.
\item It is always interesting to understand the primitive elements in a Hopf
algebra. Clearly, from the results in \cite{4TR}, there are relations between
Feynman diagrams without subdivergences. This should be a major focus
of future work.
\end{itemize}
\section*{Acknowledgements}
Foremost, I have to thank David Broadhurst for 
a long and ongoing collaboration. It is a pleasure
to thank the organizers of the
1997 UK Institute in St.Andrews for support and hospitality.
I also like to thank the students attending my recent course 
on {\em Knots and Renormalization}, especially
M.Mertens, O.Welzel, R.Kreckel, U.~Heinzmann
and V.Kleinschmidt, for lots
of interesting questions and discussions.
I also thank O.~Welzel, M.~Mertens and U.~Heinzmann for
careful reading of the manuscript, and the I.H.E.S.(Bures-sur-Yvette)
for hospitality.


\begin{thebibliography}{99}
\bibitem{Collins}
J.Collins, {\em Renormalization}, Cambridge UP (1984).
\bibitem{Zimm}
W.~Zimmermann, Comm.Math.Phys.15 (1969) 208.
\bibitem{Kassel}
C.~Kassel, {\em Quantum Groups}, Springer Verlag (1995).
\bibitem{km}
D.~Kreimer, M.~Mertens, {\em quasi Hopf Algebras from Renormalization},
in preparation.
\bibitem{mills}
R.~Mills, {\em Tutorial on Infinities in QED},
in {\em Renormalization- From Lorentz to Landau (and beyond)},
L.M.~Brown, Ed., Springer 1993.
\bibitem{habil}
D.~Kreimer, Habilschrift: {\em Renormalization 
and Knot Theory}, J.Knot Th.Ram.6 (1997) 479-581,
q-alg/9607022.
\bibitem{bdk}
D.J.~Broadhurst, R.~Delbourgo, D.~Kreimer, Phys.Lett.{\bf B366} (1996) 421.
\bibitem{Bob}
R.~Delbourgo, A.~Kalloniatis, G.~Thompson,
Phys.Rev.{\bf D54} (1996) 5373;\\
R.~Delbourgo, D.~Elliott, D.S.~McAnally,
Phys.Rev.{\bf D55} (1997) 5230.
\bibitem{db}
D.J.~Broadhurst, 
{\em On the enumeration of irreducible
k-fold Euler sums and their 
roles in knot theory and field theory},
to appear in J.Math.Phys.,
Open Univ.~preprint OUT-4102-62 (1996), hep-th/9604128.
\bibitem{plb}
D.~Kreimer, Phys.Lett.{\bf B354} (1995) 117.
\bibitem{pisa}
D.J.~Broadhurst, D.~Kreimer, Int.J.of Mod.Phys.{\bf C6} (1995) 519.
\bibitem{bgk}
D.J.~Broadhurst, J.A.~Gracey, D.~Kreimer, 
Z.Phys.{\bf C75} (1997) 559.\\
D.J.~Broadhurst, A.~Kotikov,
{\em Compact analytical form for non-zeta terms in critical
exponents at order $1/N^3$},
subm.~to Phys.Lett.B, hep-th/9612013.
\bibitem{bk15}
D.J.~Broadhurst, D.~Kreimer, 
Phys.Lett.{\bf B393} (1997) 403.
\bibitem{BBB}
J.M.~Borwein, D.A.~Bradley, D.J.~Broadhurst,
Electronic J.\ Combinatorics {\bf 4} (1997) R5.
\bibitem{db2}
D.J.~Broadhurst, {\em Conjectured enumeration of irreducible multiple
zeta values, from knots and Feynman diagrams},
subm.~to Phys.Lett.B, hep-th/9612012.
\bibitem{Zag}
D.~Zagier, in Proc.First European Congress Math.~(Birkh\"auser,
Boston, 1994) Vol.II, pp 497-512; {\em Multiple
Zeta Values}, in preparation.
\bibitem{dbn}
D.~Bar-Natan, Topology 34,2 (1995) 423;\\
{\em Non-associative Tangles}, Georgia Int.Topology Conf.Proc.
\bibitem{book}
D.~Kreimer, {\em Knots and Feynman Diagrams}
(Cambridge UP, in preparation).
\bibitem{4TR}
D.~Kreimer, {\em Weight Systems from Feynman Diagrams},
J.Knot Th.Ram.7 (1998) 61-85, hep-th/9612010.
\bibitem{bk4}
D.J.~Broadhurst, D.~Kreimer, {\em Feynman Diagrams as a weight system:
four-loop test of a four-term relation},
subm.~to Phys.Lett.B, hep-th/9612011.
\bibitem{sub}
D.~Kreimer, {\em On Knots in Subdivergent Diagrams},
to appear in Z.Phys.{\bf C}, hep-th/9610128.
\bibitem{bkhopf}
D.J.~Broadhurst, D.~Kreimer, in preparation.
\bibitem{sloane}
N.~Sloane, {\em Online Encyclopedia of Integer Sequences},  
series A027881,\\ http://www.research.att.com/$\sim$njas/sequences.
\end{thebibliography}
\end{document}